\documentclass[letterpaper,10pt]{IEEEtran}

\usepackage{amsmath,mathrsfs}
\usepackage{amssymb}
\usepackage{color}
\usepackage[all]{xy}
\usepackage{graphicx,xspace,bm}
\usepackage[caption = false]{subfig}
\usepackage{url}
\usepackage[lined,linesnumbered,algoruled]{algorithm2e}
\usepackage{setspace}
\usepackage{url}

\newtheorem{theorem}{Theorem}[section]

\newtheorem{assumption}[theorem]{Assumption}
\newtheorem{lemma}[theorem]{Lemma}

\newtheorem{remark}[theorem]{Remark}

\newtheorem{proposition}[theorem]{Proposition}  

\newcommand{\nout}{\upscr{N}{out}}

\newcommand{\Dout}{\subscr{\mathsf{D}}{out}}

\newcommand{\Adj}{\mathsf{A}}
\newcommand{\Lap}{\mathsf{L}}
\newcommand{\vertices}{\mathcal{V}}
\newcommand{\edges}{\mathcal{E}}

\newcommand{\Bgraph}{\mathcal{G}}

\renewcommand{\SS}{\mathcal{S}}
\newcommand{\KK}{\mathcal{K}}

\newcommand{\FF}{\mathcal{F}}
\newcommand{\GG}{\mathcal{G}}

\newcommand{\HH}{\mathcal{H}}

\newcommand{\OO}{\mathcal{O}}
\newcommand{\EE}{\mathcal{E}}
\newcommand{\MM}{\mathcal{M}}

\newcommand{\spn}{\mathrm{span}}

\newcommand{\ones}{\mathbf{1}}
\newcommand{\zeros}{\mathbf{0}}

\newcommand{\F}{F}

\newcommand{\eye}{\mathrm{I}}
\newcommand{\abs}[1]{\ensuremath{\left\lvert{#1}\right\rvert}}
\newcommand{\norm}[1]{\ensuremath{\| #1 \|}}
\newcommand{\inorm}[1]{\ensuremath{\| #1 \|}_{\infty}}
\newcommand{\bnorm}[1]{\ensuremath{\big \| #1 \big \|}}

\newcommand{\real}{{\mathbb{R}}}
\newcommand{\realpositive}{{\mathbb{R}}_{>0}}
\newcommand{\realnonnegative}{{\mathbb{R}}_{\ge 0}}

\newcommand{\integerspositive}{\mathbb{Z}_{\geq 1}}

\newcommand{\eps}{\epsilon}

\newcommand{\until}[1]{\{1,\dots,#1\}}
\newcommand{\map}[3]{#1:#2 \rightarrow #3}
\newcommand{\setmap}[3]{#1:#2 \rightrightarrows #3}

\newcommand{\rrarrows}{\rightrightarrows}

\newcommand{\setdef}[2]{\{#1 \; | \; #2\}}
\newcommand{\cupdot}{\mathbin{\mathaccent\cdot\cup}}

\newcommand{\SetLie}{{\mathcal{L}}}
\newcommand{\gradient}{\nabla}
\newcommand{\Eq}[1]{\mathrm{Eq}(#1)}

\newcommand{\DEDS}{{\rm DEDS}\xspace}

\newcommand{\FFDEDS}{\FF_{\mathrm{DEDS}}}
\newcommand{\FFDEDSo}{\FF_{\mathrm{DEDS}}^{*}}

\newcommand{\dac}{\texttt{dac}\xspace}

\newcommand{\dacLapgg}{X_{\texttt{dac+}(\Lap\partial,\partial)}}
\newcommand{\daclgg}{\texttt{dac+}$(\Lap\partial,\partial)$\xspace}

\newcommand{\Hglobal}{\mathfrak{M}_g}
\newcommand{\Homega}{\mathfrak{M}_o}

\newcommand{\FFa}{\FF_{\mathrm{aug}}^*}
\newcommand{\bFFa}{\overline{\FF}_{\mathrm{aug}}}

\newcommand{\lm}{\lambda}

\DeclareMathAlphabet{\mathpzc}{OT1}{pzc}{m}{it}

\newcommand\subscr[2]{#1_{\textup{#2}}}
\newcommand\upscr[2]{#1^{\textup{#2}}}

\newcommand{\oprocendsymbol}{\hbox{$\bullet$}}
\newcommand{\oprocend}{\relax\ifmmode\else\unskip\hfill\fi\oprocendsymbol}

\newcommand{\longthmtitle}[1]{\mbox{}\textup{\textsl{(#1):}}}

\newcommand{\horizon}{\mathfrak{h}}

\renewcommand{\theenumi}{(\roman{enumi})}

\newcommand{\myclearpage}{\clearpage}
\renewcommand{\myclearpage}{}

\begin{document}

\title{Distributed coordination of DERs with storage
  \\
  for dynamic economic dispatch\thanks{A preliminary version appeared
    as~\cite{AC-JC:15-cdc} at the IEEE Conference on Decision and
    Control.}}

\author{Ashish Cherukuri \qquad Jorge Cort\'{e}s\thanks{Ashish
    Cherukuri and Jorge~Cort\'{e}s are with the Department of
    Mechanical and Aerospace Engineering, University of California,
    San Diego, \texttt{\{acheruku,cortes\}@ucsd.edu}.}}

\maketitle

\begin{abstract}
  This paper considers the dynamic economic dispatch problem for a
  group of distributed energy resources (DERs) with storage that communicate over a
  weight-balanced strongly connected digraph.  The objective is to
  collectively meet a certain load profile over a finite time horizon
  while minimizing the aggregate cost.  At each time slot, each
  DER decides on the amount of generated power, the amount sent
  to/drawn from the storage unit, and the amount injected into the
  grid to satisfy the load.  Additional constraints include bounds on
  the amount of generated power, ramp constraints on the difference in
  generation across successive time slots, and bounds on the amount of
  power in storage.  We synthesize a provably-correct distributed
  algorithm that solves the resulting finite-horizon optimization
  problem starting from any initial condition. Our design consists of
  two interconnected systems, one estimating the mismatch
  between the injection and the total load at each time slot, and
  another using this estimate to reduce the mismatch and optimize the
  total cost of generation while meeting the
  constraints. 
\end{abstract}

\section{Introduction}\label{sec:Intro}

The current electricity grid is up for a major transformation to
enable the widespread integration of distributed energy resources and
flexible loads to improve efficiency and reduce emissions without
affecting reliability and performance. This presents the need for
novel coordinated control and optimization strategies which, along
with suitable architectures, can handle uncertainties and variability,
are fault-tolerant and robust, and preserve privacy. With this context
in mind, our objective here is to provide a distributed algorithmic
solution to the dynamic economic dispatch problem with storage.  We
see the availability of such strategies as a necessary building block
in realizing the vision of the future grid.

\emph{Literature review:} Static economic dispatch (SED) involves a
group of generators collectively meeting a specified load for a single
time slot while minimizing the total cost and respecting individual
constraints. In recent years, distributed generation has motivated the
shift from traditional solutions of the SED problem to decentralized
ones, see e.g.,~\cite{ADDG-STC-CNH:12,SK-GH:12,WZ-WL-XW-LL-FF:14} and
our own
work~\cite{AC-JC:15-tcns,AC-JC:14-auto}.  As argued
in~\cite{XX-AME:10,MDI-LX-JJ:11}, the dynamic version of the problem,
termed dynamic economic dispatch (DED), results in better grid control
as it optimally plans generation across a time horizon, specifically
taking into account ramp limits and variability of power commitment
from renewable sources.  Conventional solution methods to the DED
problem are centralized~\cite{XX-AME:10}.  Recent
works~\cite{MDI-LX-JJ:11,XX-JZ-AE:11} have employed model predictive
control (MPC)-based algorithms to deal more effectively with complex
constraints and uncertainty, but the resulting methods are still
centralized and do not provide theoretical guarantees on the
optimality of the solution.  The work~\cite{ZL-WW-BZ-HS-QG:13}
proposes a Lagrangian relaxation method to solve the DED problem, but
the implementation requires a master agent that communicates with and
coordinates the generators. MPC methods have also been employed
by~\cite{AH-JM-HM-HD:13} in the dynamic economic dispatch with storage
(\DEDS) problem, which adds storage units to the DED problem to lower
the total cost, meet uncertain demand under uncertain generation, and
smooth out the generation profile across time.  The stochastic version
of the \DEDS problem adds uncertainty in demand and generation by
renewables. Algorithmic solutions for this problem put the emphasis on
breaking down the complexity to speed up convergence for large-scale
problems and include stochastic MPC~\cite{DZ-GH:14}, dual
decomposition~\cite{YZ-NG-GBG:13}, and optimal condition
decomposition~\cite{AS-AR-AK:15} methods. However, these methods are
either centralized or need a coordinating central master.

\emph{Statement of contributions:} Our starting point is the
formulation of the \DEDS problem for a group of power DERs 
communicating over a weight-balanced strongly connected digraph. Since
the cost functions are convex and all constraints are linear, the
problem is convex in its decision variables, which are the power to be
injected and the power to be sent to storage by each DER at each
time slot. Using exact penalty functions, we reformulate the \DEDS
problem as an equivalent optimization that retains equality
constraints but removes inequality ones.  The structure of the
modified problem guides our design of the provably-correct distributed
strategy termed ``dynamic average consensus (\dac) + Laplacian
nonsmooth gradient ($\Lap \partial$) + nonsmooth gradient
($\partial$)'' dynamics to solve the \DEDS problem starting from any
initial condition.  This algorithm consists of two interconnected
systems. A first block allows DERs to track, using \dac, the
mismatch between the current total power injected and the load for
each time slot of the planning horizon. A second block has two
components, one that minimizes the total cost while keeping the total
injection constant (employing Laplacian-nonsmooth-gradient dynamics on
injection variables and nonsmooth-gradient dynamics on storage
variables) and an error-correcting component that uses the mismatch
signal estimated by the first block to adjust, exponentially fast, the
total injection towards the load for each time slot.

\emph{Notation:} Let $\real$, $\realnonnegative$, $\realpositive$,
$\integerspositive$ denote the set of real, nonnegative real, positive
real, and positive integer numbers, respectively.  The $2$- and
$\infty$-norm on $\real^n$ are denoted by $\norm{\cdot}$ and
$\norm{\cdot}_{\infty}$, respectively.  We let $B(x,\delta)$
denote the open ball centered at $x \in \real^n$ with radius $\delta >
0$. Given $r \in \real$, we denote $\HH_r = \setdef{x \in
  \real^n}{\ones_n^\top x = r}$.
For a symmetric matrix $A \in \real^{n \times n}$, the minimum and
maximum eigenvalues of $A$ are $\lambda_{\min}(A)$ and
$\lambda_{\max}(A)$. The Kronecker product of $A \in \real^{n \times
  m}$ and $B \in \real^{p \times q}$ is $A \otimes B \in \real^{np
  \times mq}$.  We use $\zeros_n = (0,\ldots,0) \in \real^n$,
$\ones_n=(1,\ldots,1) \in \real^n$, and $\eye_n \in \real^{n \times
  n}$ for the identity matrix.  For $x\in \real^n$ and $y \in
\real^m$, the vector $(x;y) \in \real^{n+m}$ denotes the
concatenation.  Given $x,y\in \real^n$, $x_i$ denotes the $i$-th
component of $x$, and $x \le y$ denotes $x_i \le y_i$ for $i \in
\until{n}$.  For $\horizon > 0$, given $y \in \real^{n\horizon}$ and
$k \in \until{\horizon}$, the vector containing the $nk-n+1$ to $nk$
components of $y$ is $y^{(k)} \in \real^n$, and so, $y =
(y^{(1)};y^{(2)}; \dots; y^{(\horizon)})$. We let $[u]^{+} = \max
\{0,u\}$ for~$u \in \real$.  A set-valued map
$\setmap{f}{\real^{n}}{\real^{m}}$ associates to each point in
$\real^{n}$ a set in $\real^{m}$.

\myclearpage
\section{Preliminaries}\label{se:Prelim}

This section introduces concepts from graph theory, nonsmooth
analysis, differential inclusions, and optimization.

\vspace*{0.5ex}%
\emph{Graph theory:} 
Following~\cite{FB-JC-SM:08cor}, a \emph{weighted directed graph}, is
a triplet $\Bgraph=(\vertices,\edges,\Adj)$, where $\vertices$ is the
vertex set, $ \edges \subseteq \vertices\times \vertices $ is the edge
set, and $ \Adj \in \mathbb{R}^{n\times n}_{\geq0} $ is the
\emph{adjacency matrix} with the property that $a_{ij}>0 $ if $
(v_i,v_j)\in \edges $ and $ a_{ij}=0 $, otherwise.  A path is an
ordered sequence of vertices such that any consecutive pair of
vertices is an edge.  A digraph is \emph{strongly connected} if there
is a path between any pair of distinct vertices.  For a vertex $v_i$,
$\nout(v_i) = \setdef{v_j \in \vertices}{(v_i, v_j) \in \edges}$ is
the set of its out-neighbors.
The \emph{Laplacian} matrix is $ \Lap = \Dout -\Adj$, where $\Dout$ is
the diagonal matrix defined by $(\Dout)_{ii}= \sum_{j=1}^{n}a_{ij} $,
for all $ i \in \{1,\ldots,n\}$. Note that $\Lap\ones_n=0 $.  If~$
\Bgraph $ is strongly connected, then zero is a simple eigenvalue of
$\Lap$.  $\Bgraph$ is \emph{weight-balanced}
iff $ \ones_n^\top \Lap =0 $ iff $ \Lap+\Lap^\top$ is positive
semidefinite.  If~$\Bgraph $ is weight-balanced and strongly
connected, then zero is a simple eigenvalue of $\Lap+\Lap^\top $ and,
for $x \in \real^n$,
\begin{equation}\label{eq:LapBound}
  \lambda_2(\Lap + \Lap^\top) \bnorm{x-\frac{1}{n}(\ones_n^\top x)\ones_n}^2
  \le x^\top (\Lap + \Lap^\top)x,
\end{equation}
with $\lambda_2(\Lap + \Lap^\top)$ the smallest non-zero eigenvalue of
$\Lap + \Lap^\top$.

\vspace*{.5ex}%
\emph{Nonsmooth analysis:} 
Here, we introduce some notions on nonsmooth analysis
from~\cite{JC:08-csm-yo}.  A function $f:\real^n \rightarrow
\real^m$ is \emph{locally Lipschitz} at $x \in \real^n$ if there exist
$L, \epsilon \in \realpositive$ such that
  $\norm{f(y) - f(y')} \le L\norm{y - y'}$,
for all $y, y'\in B(x,\epsilon)$.
A function $f:\real^n \rightarrow \real$ is \emph{regular} at $x \in
\real^n$ if, for all $v \in \real^n$, the right directional derivative
and the generalized directional derivative of $f$ at $x$ along the
direction $v$ coincide,
see~\cite{JC:08-csm-yo} for these definitions. 
A convex function is regular.  A set-valued map $\HH:\real^n \rrarrows
\real^n$ is \emph{upper semicontinuous} at $x \in \real^n$ if, for all
$\epsilon \in \realpositive$, there exists $\delta \in \realpositive$
such that $\HH(y) \subset \HH(x) + B(0,\epsilon)$ for all $y \in
B(x,\delta)$. Also, $\HH$ is \emph{locally bounded} at $x \in \real^n$
if there exist $\epsilon, \delta \in \realpositive$ such that
$\norm{z} \le \epsilon$ for all $z \in \HH(y)$, and all $y \in
B(x,\delta)$.  Given a locally Lipschitz function $f:\real^n
\rightarrow \real$, let $\Omega_f$ be the set (of measure zero) of
points where $f$ is not differentiable. The \emph{generalized
  gradient} $\partial f: \real^n \rrarrows \real^n$
of~$f$~is 
\begin{equation*}
  \partial f(x) = \mathrm{co} \setdef{ \lim_{i \rightarrow \infty} 
    \gradient f(x_i)}{ x_i \rightarrow x, x_i \notin S \cup \Omega_f},
\end{equation*}
where $\mathrm{co}$ is the convex hull and $S \subset \real^n$ is any
set of measure zero.  The set-valued map $\partial f$ is locally
bounded, upper semicontinuous, and takes non-empty, compact, and
convex values. 
For a function $\map{f}{\real^n \times
\real^m} {\real}$,  $(x,y) \mapsto f(x,y)$,
the partial generalized gradient with respect to $x$ and $y$ are
denoted by $\partial_x f$ and $\partial_y f$, respectively. 

\vspace*{.5ex}%
\emph{Differential inclusions:} 
We gather here tools from~\cite{JC:08-csm-yo,AC-JC:14-auto} to analyze
the stability properties of differential inclusions,
\begin{equation}\label{eq:ddsys}
  \dot x \in \F(x) ,
\end{equation}
where $\F: \real^n \rrarrows \real^n$ is a set-valued map.  A solution
of~\eqref{eq:ddsys} on $[0,T] \subset \real$ is an absolutely
continuous map $x:[0,T]\rightarrow \real^n$ that
satisfies~\eqref{eq:ddsys} for almost all $t \in [0,T]$.  If the
set-valued map $\F$ is locally bounded, upper semicontinuous, and
takes non-empty, compact, and convex values, then the existence of
solutions is guaranteed. The set of equilibria of~\eqref{eq:ddsys} is
$\Eq{\F} = \setdef{x \in \real^n}{0 \in \F(x) }$.
Given a locally Lipschitz function $W: \real^n \rightarrow \real$, the
\emph{set-valued Lie derivative} $\SetLie_{\F}W: \real^n \rrarrows
\real$ of $W$ with respect to~\eqref{eq:ddsys} at $x \in \real^n$ is
\begin{align*}
  \SetLie_{\F}W = \setdef{a \in \real}{\exists v \in
    \F(x)  \text{ s.t. }  \zeta^\top v=a, \, 
    \forall 
    \zeta \in \partial W(x)} .
\end{align*}
The \emph{$\omega$-limit set} of a trajectory $t \mapsto \varphi(t)$,
$\varphi(0) \in \real^n$ of~\eqref{eq:ddsys}, denoted
$\Omega(\varphi)$, is the set of all points $y \in \real^n$ for which
there exists a sequence of times $\{t_k\}_{k=1}^\infty$ with $t_k \to
\infty$ such that $\lim_{k \to \infty} \varphi(t_k) = y$. If the
trajectory is bounded, then the $\omega$-limit set is nonempty,
compact, connected.  The next result from~\cite{AC-JC:14-auto} is a
refinement of the LaSalle Invariance Principle for differential
inclusions that establishes convergence
of~\eqref{eq:ddsys}.

\begin{proposition}\longthmtitle{Refined LaSalle Invariance Principle
    for differential inclusions}\label{pr:refined-lasalle-nonsmooth}
  Let $\setmap{\F}{\real^n}{\real^n}$ be upper semicontinuous, taking
  nonempty, convex, and compact values everywhere in $\real^n$.
  Let $t \mapsto \varphi(t)$ be a bounded solution of~\eqref{eq:ddsys}
  whose $\omega$-limit set $\Omega(\varphi)$ is contained in $\SS
  \subset \real^n$, a closed embedded submanifold of $\real^n$.  Let
  $\OO$ be an open neighborhood of $\SS$ where a locally Lipschitz,
  regular function $\map{W}{\OO}{\real}$ is defined.  Then,
  $\Omega(\varphi) \subset \EE$ if the following holds,
  \begin{enumerate}
  \item $\!\EE \!=\! \setdef{x \in \SS\!}{\!0 \in \SetLie_{\F} W(x) }$ belongs
    to a level set of~$W$
    \label{as:refined-lasalle-1}
  \item for any compact set $\MM \subset \SS$ with $\MM \cap \EE =
    \emptyset$, there exists a compact neighborhood $\MM_c$ of $\MM$
    in $\real^n$ and $\delta < 0$ such that $\sup_{x \in \MM_c} \max
    \SetLie_{\F} W(x) \le \delta$.
    \label{as:refined-lasalle-2}
  \end{enumerate}
\end{proposition}

\vspace*{.5ex}%
\emph{Constrained optimization and exact penalty
  functions:} 
  Here, we introduce some notions on constrained convex optimization
following~\cite{SB-LV:09,DPB:75b}.  Consider the optimization problem,
\begin{subequations}\label{eq:GenConsOpt}
  \begin{align}
    \mathrm{minimize} \quad & f(x), \label{eq:GenConsObjective}
    \\
    \text{subject to} \quad & g(x) \le \zeros_m, \quad h(x)
    =\zeros_p, \label{eq:GenConsConstraints}
  \end{align}
\end{subequations}
where $f:\real^n \rightarrow \real$, $g:\real^n \rightarrow \real^m$,
are continuously differentiable and convex, 
and $h:\real^n \rightarrow \real^p$ with $p\le n$ is affine.  
The \emph{refined Slater condition} is satisfied
by~\eqref{eq:GenConsOpt} if there exists $x \in \real^n$ such that
$h(x) = \zeros_p$, $g(x) \le \zeros_m$, and $g_i(x) < 0$ for all
nonaffine functions $g_i$. The refined Slater condition implies that
strong duality holds.
A point $x \in \real^n$ is a Karush-Kuhn-Tucker (KKT) point
of~\eqref{eq:GenConsOpt} if there exist Lagrange multipliers $\lambda
\in \realnonnegative^m$ and $\nu \in \real^p$ such that
\begin{align*}
  & g(x) \le \zeros_m, \quad h(x) = \zeros_p, \quad \lambda^\top g(x) =
  0,
  \\
  & \gradient f(x)+ \sum_{i=1}^m \lambda_i \gradient g_i(x) +
  \sum_{i=1}^p \nu_i \gradient h_i(x)  = 0.
\end{align*}
If strong duality holds then, a point is a solution
of~\eqref{eq:GenConsOpt} iff it is a KKT point.  The
optimization~\eqref{eq:GenConsOpt} satisfies the \emph{strong Slater
  condition} with parameter $\rho \in \realpositive$ and feasible
point $x^\rho \in \real^n$ if $g(x^\rho) < - \rho \ones_m$ and
$h(x^\rho) = \zeros_p$.

\begin{lemma}\longthmtitle{Bound on Lagrange
    multiplier~\cite[Remark 2.3.3]{JBHU-CL:93}}\label{le:lagrange-bound}
  If~\eqref{eq:GenConsOpt} satisfies the strong Slater condition with
  parameter $\rho \in \realpositive$ and feasible point $x^\rho \in
  \real^n$, then any primal-dual optimizer $(x,\lm,\nu)$
  of~\eqref{eq:GenConsOpt} satisfies
  \begin{align*}
    \inorm{\lm} \le \frac{f(x^\rho) - f(x)}{\rho}.
  \end{align*}
\end{lemma}

We are interested in 
eliminating the inequality constraints in~\eqref{eq:GenConsOpt} while
keeping the equality constraints intact.
To this end, we use~\cite{DPB:75b} to construct a nonsmooth exact
penalty function $f^{\epsilon}: \real^n \rightarrow \real$, given as $
f^{\epsilon}(x) = f(x) + \frac{1}{\epsilon} \sum_{i=1}^m [g_i(x)]^+$,
with $\epsilon >0$, and define the minimization problem
\begin{subequations}\label{eq:ExactPenalty}
  \begin{align}
    \mathrm{minimize} \quad & f^{\epsilon}(x),
  \label{eq:ExactPenalty1} 
  \\
  \text{subject to} \quad & h(x) = \zeros_p. \label{eq:ExactPenalty2}
  \end{align}
\end{subequations}
Note that $f^\eps$ is convex as $f$ and $t \mapsto \frac{1}{\epsilon}
[t]^+$ are convex. Hence, 
the problem~\eqref{eq:ExactPenalty} is convex.
The following
result, see e.g.~\cite[Proposition 1]{DPB:75b}, identifies conditions
under which the solutions of the problems~\eqref{eq:GenConsOpt}
and~\eqref{eq:ExactPenalty} coincide.

\begin{proposition}\longthmtitle{Equivalence of~\eqref{eq:GenConsOpt}
    and~\eqref{eq:ExactPenalty}}\label{pr:EquivalenceExactPenalty}
  Assume~\eqref{eq:GenConsOpt} has nonempty, compact solution set, and
  satisfies the refined Slater condition. Then,~\eqref{eq:GenConsOpt}
  and~\eqref{eq:ExactPenalty} have the same solutions if $
  \frac{1}{\epsilon} > \norm{\lambda}_\infty$, for some Lagrange
  multiplier $\lambda \in \realnonnegative^m$
  of~\eqref{eq:GenConsOpt}.
\end{proposition}

\myclearpage
\section{Problem statement}\label{sec:problem}

Consider a network of $n \in \integerspositive$ distributed energy
resources (DERs) whose communication topology is a strongly connected
and weight-balanced digraph $\Bgraph=(\vertices,\edges,\Adj)$.  For
simplicity, we assume DERs to be generator units.  In our discussion,
DERs can also be flexible loads (where the cost function corresponds
to the negative of the load utility function).
An edge $(i,j)$ represents the capability of unit $j$ to transmit
information to unit~$i$. Each unit $i$ is equipped with storage
capabilities with minimum~$C^m_i \in \realnonnegative$ and maximum
$C^M_i \in \realpositive$ capacities.  The network collectively aims
to meet a power demand profile during a finite-time horizon $\KK =
\until{\horizon}$ specified by $l \in \realpositive^\horizon$, that
is, $l^{(k)}$ is the demand at time slot $k \in \KK$.  This demand can
either correspond to a load requested from an outside entity, denoted
$L^{(k)} \ge 0$ for slot $k$, or each DER $i$ might have to
satisfy a load at the bus it is connected to, denoted
$\tilde{l}_i^{(k)} \ge 0$ for slot $k$. Thus, for each $k \in \KK$,
$l^{(k)} = L^{(k)} + \sum_{i=1}^n \tilde{l}_i^{(k)}$. We assume that
the external demand $L = (L^{(1)}, \dots, L^{(\horizon)}) \in
\realnonnegative^{\horizon}$ is known to an arbitrarily selected
unit~$r \in \until{n}$, whereas the demand at bus $i$, $\tilde{l}_i =
(\tilde{l}_i^{(1)}, \dots, \tilde{l}_i^{(\horizon)}) \in
\realnonnegative^{\horizon}$, is known to unit~$i$.  For convenience,
$\tilde{l} = (\tilde{l}^{(1)}, \dots, \tilde{l}^{\horizon})$, where
$\tilde{l}^{(k)} = (\tilde{l}^{(k)}_1, \dots, \tilde{l}^{(k)}_n)$
collects the load known to each unit at slot~$k \in \KK$.
Along with load satisfaction, the group also
aims to minimize the total cost of generation and to satisfy the
individual physical constraints for each DER. We make these
elements precise next.

Each unit $i$ decides at every time slot $k$ in $\KK$ the amount of
power it generates, the portion $I^{(k)}_i \in \real$ of it that it
injects into the grid to meet the load, and the remaining part
$S^{(k)}_i \in \real$ that it sends to the storage unit.  The power
generated by $i$ at $k$ is then $I^{(k)}_i + S^{(k)}_i$.  We denote by
$I^{(k)} = (I_1^{(k)}, \dots, I_n^{(k)}) \in \real^n$ and $S^{(k)} =
(S_1^{(k)}, \dots, S_n^{(k)}) \in \real^n$ the collective injected and
stored power at time $k$, respectively.  The load satisfaction is then
expressed as $\ones_n^\top I^{(k)} = l^{(k)} = L^{(k)} + \ones_n^\top
\tilde{l}^{(k)}$, for all $k \in \KK$.  The cost
$f_i^{(k)}(I_i^{(k)}+S_i^{(k)})$ of power generation $I^{(k)}_i +
S^{(k)}_i$ by unit $i$ at time $k$ is specified by the function
$f_i^{(k)}:\real \rightarrow \realnonnegative$, which we assume convex
and continuously
differentiable.  
Given $(I^{(k)},S^{(k)})$, the cost incurred by the network at time
slot $k$ is
\begin{align*}
  f^{(k)}(I^{(k)}+S^{(k)}) =
  \sum_{i=1}^{n}f_i^{(k)}(I_i^{(k)}+S_i^{(k)}) .
\end{align*}
The cumulative cost of generation for the network across the time
horizon is $\map{f}{\real^{n\horizon}}{\realnonnegative}$, $f(x) =
\sum_{k=1}^\horizon f^{(k)}(x^{(k)})$.  Given injection $I = (I^{(1)},
\dots, I^{(\horizon)}) \in \real^{n\horizon}$ and storage $S =
(S^{(1)}, \dots, S^{(\horizon)}) \in \real^{n\horizon}$ values, the
total network cost is 
\begin{align*}
  f(I+S) = \sum_{k=1}^\horizon f^{(k)} (I^{(k)}+S^{(k)}).
\end{align*}
The functions $\{f^{(k)}\}_{k \in \KK}$ and $f$ are also convex and
continuously differentiable.  Next, we describe the physical
constraints on the DERs. Each unit's power must belong to
the range~$[P^m_i, P^M_i] \subset \realpositive$, representing lower
and upper bounds on the amount of power it can generate at each time
slot.  Each unit $i$ also respects upper and lower ramp
constraints: the change in the generation level from any time slot $k$
to $k+1$ is upper and lower bounded by $R^u_i$ and $-R^l_i$,
respectively, with $ R^u_i$, $R^l_i \in \realpositive $.  At each time
slot, the power injected into the grid by each unit must be
nonnegative, i.e., $I^{(k)}_i \ge 0$.  Furthermore, the amount of
power stored in any storage unit $i$ at any time slot $k \in \KK$ must
belong to the range~$[C^m_i, C^M_i]$.  Finally, we assume that at the
beginning of the time slot $k=1$, each storage unit $i$ starts with
some stored power $S^{(0)}_i \in [C^m_i,C^M_i]$.  With the above
model, the \emph{dynamic economic dispatch with storage} (\DEDS)
problem is formally defined by the following convex optimization
problem,
\begin{subequations}\label{eq:deds}
  \begin{align}
    \underset{(I,S) \in \real^{2n\horizon}}{\text{minimize}} & \quad 
    f(I+S), \label{eq:conobjective-s}
    \\
    \text{subject to}&  \, \, \,  \text{for } k \in \KK, \notag
    \\
    & \quad \ones_n^\top I^{(k)} = l^{(k)}, \label{eq:load-cond}
    \\
    & \quad P^{m} \le I^{(k)} + S^{(k)} \le P^{M}, \label{eq:box-cons}
    \\
    & \textstyle \quad C^m \le S^{(0)} + \sum_{k' = 1}^k S^{(k')} \le C^M,
    \label{eq:storage-cons}
    \\
    & \quad \zeros_n \le I^{(k)}, \label{eq:injection-cons}
    \\
    & \, \, \, \text{for } k \in \KK \setminus \{\horizon\}, \notag
    \\
    & \!-\! R^l \le I^{(k+1)} \!+\! S^{(k+1)} \!-\! I^{(k)} \!-\! S^{(k)}
    \le R^u. \label{eq:ramp-cons}
  \end{align}
\end{subequations}
We refer to~\eqref{eq:load-cond}--\eqref{eq:ramp-cons} as the
\emph{load conditions}, \emph{box constraints}, \emph{storage limits},
\emph{injection constraints}, and \emph{ramp constraints},
respectively.  We denote by $\FFDEDS$ and $\FFDEDSo$ the feasibility
set and the solution set of the \DEDS problem~\eqref{eq:deds},
respectively, and assume them to be nonempty.  Since $\FFDEDS$ is
compact, so is $\FFDEDSo$.  Moreover, the refined Slater condition is
satisfied for \DEDS as all the
constraints~\eqref{eq:load-cond}--\eqref{eq:ramp-cons} are affine in
the decision variables. Additionally, we assume that the \DEDS problem
satisfies the strong Slater condition with parameter $\rho \in
\realpositive$ and feasible point $(I^\rho,S^\rho) \in
\real^{2n\horizon}$.

\begin{remark}\longthmtitle{General setup for
    storage}\label{re:storage-subset}
  The \DEDS formulation above can be modified to consider scenarios
  where only some DERs $\vertices_{gs}$ are equipped with
  storage  and others $\vertices_{g}$ are not, with
  $\until{n} = \vertices_{gs} \cupdot \vertices_{g}$.  The formulation
  can also be extended to consider the cost of storage,
  inefficiencies, and constraints on (dis)charging of the storage
  units, as in~\cite{AH-JM-HM-HD:13,YZ-NG-GBG:13}.  These factors
  either affect the constraint~\eqref{eq:storage-cons}, add additional
  conditions on the storage variables, or modify the objective
  function.  As long as the resulting cost and constraints are convex
  in $S$, all these can be treated within~\eqref{eq:deds} without
  affecting the design methodology.
  \oprocend
\end{remark}

Our aim is to design a distributed algorithm that allows the network
interacting over~$\GG$ to solve the \DEDS problem.

\myclearpage
\section{Distributed algorithmic solution}
\label{sec:main-result}

We describe here the distributed algorithm that asymptotically finds
the optimizers of the \DEDS problem. Our design strategy builds on an
alternative formulation of the optimization problem using penalty
functions (cf. Section~\ref{sec:main-result}-A). This allows us to get
rid of the inequality constraints, resulting into an optimization
whose structure guides our algorithmic design
(cf. Section~\ref{sec:main-result}-B).

\vspace*{0.5ex}
\emph{A. Alternative formulation of
the \DEDS problem:} 
The procedure here follows closely the theory of exact penalty
functions outlined in Section~\ref{se:Prelim}.
For an $\epsilon \in \realpositive$, consider the modified cost
function $\map{f^{\eps}}{\real^{n\horizon} \times
\real^{n\horizon}}{\realnonnegative}$,
\begin{align*}
  f^{\eps}(I,& S)  = f(I+S) + \frac{1}{\eps} \Bigl( \sum_{k = 1}^{\horizon} 
  \ones_n^\top \bigl( [T^{(k)}_1]^+ + [T^{(k)}_2]^+ + [T^{(k)}_3]^+ 
  \\
  &  +
  [T^{(k)}_4]^+ + [T^{(k)}_5]^+ \bigr) 
   + \sum_{k=1}^{\horizon-1} \ones_n^\top \bigl( [T^{(k)}_6]^+ +
  [T^{(k)}_7]^+ \bigr) \Bigr),
\end{align*}
where
\begin{align}\label{eq:T-defs}
  &  T^{(k)}_1  = P^m - I^{(k)} - S^{(k)},
  \, T^{(k)}_2  = I^{(k)} + S^{(k)} - P^M, \notag
  \\
  & \textstyle T^{(k)}_3  = C^m - S^{(0)}- \sum_{k'=1}^k S^{(k')}, \notag
  \\
  & \textstyle T^{(k)}_4 = S^{(0)} + \sum_{k'=1}^k S^{(k')} - C^M,
  \, \, T^{(k)}_5  = - I^{(k)},  \notag 
  \\
  &T^{(k)}_6  = -R^l - I^{(k+1)}- S^{(k+1)}  + I^{(k)}+ S^{(k)}, \notag
  \\
  &T^{(k)}_7  = I^{(k+1)} + S^{(k+1)} - I^{(k)}- S^{(k)}  - R^u. 
\end{align}
This cost contains the penalty terms for all the inequality
constraints of the \DEDS problem.  Note that $f^\eps$ is locally
Lipschitz, jointly convex in $I$ and $S$, and regular.  Thus, the
partial generalized gradients $\partial_I f^\eps$ and $\partial_S
f^\eps$ take nonempty, convex, compact values and are locally bounded
and upper semicontinuous.  Consider the modified \DEDS problem
\begin{subequations}\label{eq:deds-m}
  \begin{align}
    \mathrm{minimize} \quad &
    f^{\eps}(I,S), \label{eq:conobjective-ms}
    \\
    \text{subject to} \quad & \ones_n^\top I^{(k)} = l^{(k)}, \,
    \forall k \in \KK. \label{eq:equalitycons-ms}
  \end{align}
\end{subequations}
The next result provides a criteria for selecting $\eps$ such that the
modified \DEDS and the \DEDS problems have the exact same
solutions. The proof is a direct application of
Lemmas~\ref{le:lagrange-bound} and~\ref{pr:EquivalenceExactPenalty}
using that the \DEDS problem satisfies the strong Slater condition
with parameter $\rho$ and feasible point~$(I^\rho,S^\rho)$.

\begin{lemma}\longthmtitle{Equivalence of \DEDS and modified \DEDS
    problems}\label{le:equiv-ded}
  Let $(I^*,S^*) \in \FFDEDSo$. Then, the optimizers of the
  problems~\eqref{eq:deds} and~\eqref{eq:deds-m} are the same for $\eps
  \in \realpositive$ satisfying
  \begin{align}\label{eq:eps-bound}
    \eps < \frac{\rho}{f(I^\rho+S^\rho) - f(I^*+S^*)}.
   \end{align}
\end{lemma}
As a consequence, if $\eps$ satisfies~\eqref{eq:eps-bound} then,
writing the Lagrangian and the KKT conditions for~\eqref{eq:deds-m} 
gives the following characterization of the
solution set of the \DEDS problem
\begin{align}\label{eq:solution-set}
  \FFDEDSo = & \setdef{(I,S) \in \real^{2n\horizon}}{\ones_n^\top I^{(k)} =
  l^{(k)} \text{ for all } k \in \KK , \notag
  \\
  & 0 \in \partial_S f^\eps(I,S), \text{ and } \exists \nu \in \real^\horizon
  \text{ such that }\notag 
  \\
  & (\nu^{(1)} \ones_n; \dots ; \nu^{(\horizon)} \ones_n) \in \partial_I
    f^\eps (I,S)}.
\end{align}
Recall that $\FFDEDSo$ is bounded.  Next, we stipulate a mild
regularity assumption on this set which implies that perturbing it by
a small parameter does not result into an unbounded set. This property
is of use in our convergence analysis later.

\begin{assumption}\longthmtitle{Regularity of $\FFDEDSo$}
  \label{as:regular} 
  For $p \in \realnonnegative$, define the map $p \mapsto \FF(p)
  \subset \real^{2n\horizon}$ as  
  \begin{align*}
    \FF(p)  = & \setdef{(I,S) \in \real^{2n\horizon}}{\abs{\ones_n^\top I^{(k)} -
    l^{(k)}} \le p \text{ for all } k \in \KK, 
    \\
    & 0 \in \partial_S f^\eps(I,S) + p B(0,1), \text{ and } \exists \nu
    \in \real^\horizon \text{ such that } 
    \\
    &  (\nu^{(1)} \ones_n; \dots ; \nu^{(\horizon)} \ones_n) \in \partial_I
    f^\eps (I,S) + p B(0,1)}.
  \end{align*}
  Note that $\FF(0) = \FFDEDSo$.
  Then, there exists a $\bar{p} > 0$ such that $\FF(p)$ is bounded for
  all $p \in [0,\bar{p})$. \oprocend
\end{assumption}

We end this section by stating a property of the generalized gradient
of $f^\eps$ that will be employed later in the analysis.

\begin{lemma}\longthmtitle{Uniform bound on the difference between
   $\partial_I f^\eps$ and $\partial_S
   f^\eps$}\label{le:size-gen-gradient}
 For $(I,S) \in \real^{2n\horizon}$, any two elements $\zeta_1 \in
 \partial_I f^\eps(I,S)$ and $\zeta_2 \in \partial_S f^\eps(I,S)$
 satisfy
 \begin{align*}
   \inorm{\zeta_1 - \zeta_2} \le \frac{\horizon+4}{\eps}.
 \end{align*}
\end{lemma}
\begin{IEEEproof}
  Write $f^\eps(I,S) = f_a(I+S) + f_b(I) + f_c(S)$ where the functions
  $f_a,f_b,f_c: \real^{n\horizon} \to \realnonnegative$ are 
  \begin{align*}
    f_a(I+S) & = f(I+S) + \frac{1}{\eps} \Bigl(\sum_{k=1}^\horizon
    \ones_n^\top ( [T_1^{(k)}]^+ + [T_2^{(k)}]^+ )
    \\
    & \qquad \qquad \quad + \sum_{k=1}^{\horizon-1} \ones_n^\top (
    [T_6^{(k)}]^+ + [T_7^{(k)}]^+) \Bigr),
    \\
    f_b(I) & = \frac{1}{\eps} \sum_{k=1}^\horizon \ones_n^\top
    [T_5^{(k)}]^+,
    \\
    f_c(S) & = \frac{1}{\eps} \sum_{k=1}^\horizon \ones_n^\top
    ([T_3^{(k)}]^+ + [T_4^{(k)}]^+).
  \end{align*}
  From the sum rule of generalized gradients~\cite{JC:08-csm-yo}, any
  element $\zeta_1 \in \partial_I f^\eps(I,S)$ can be expressed as a
  sum of the vectors $\zeta_{1,a}$ and $\zeta_{1,b} \in \real^{n\horizon}$ 
  such that $\zeta_{1,a} \in
  \partial f_a(I+S)$ and $\zeta_{1,b} \in \partial
  f_b(I)$. Similarly, $\zeta_2 = \zeta_{2,a} + \zeta_{2,c}$ where
  $\zeta_{2,a} \in \partial f_a(I+S)$ and $\zeta_{2,c} \in
  \partial f_c(S)$. By the definition of $f_b$, we get
  $\inorm{\zeta_{1,b}} \le \frac{1}{\eps}$. For the function $f_c$,
  note that for any $i \in \until{n}$ and any $k \in \KK$, either
  $([T_3^{(k)}]^+)_i$ is zero or $([T_4^{(k)}]^+)_i$ is zero.
  Considering extreme case, if for a particular $i$, either
  $([T_3^{(k)}]^+)_i > 0$ or $([T_4^{(k)}]^+)_i > 0$ for all $k \in
  \KK$ then, we obtain $\abs{(\zeta_{2,c})_i^{(1)}} =
  \frac{\horizon}{\eps}$. This implies that $\inorm{\zeta_{2,c}} \le
  \frac{\horizon}{\eps}$.  Now consider any two elements $\zeta_{1,a},
  \zeta_{2,a} \in \partial f_a(I+S)$. Note that for any $i \in
  \until{n}$, either $([T_1^{(k)}]^+)_i$ is zero or
  $([T_2^{(k)}]^+)_i$ is zero. Similarly, either $([T_6^{(k)}]^+)_i$
  or $([T_7^{(k)}]^+)_i$ is zero. Further, note that $I^{(k)}_i +
  S^{(k)}_i$ appears in $([T^{(k)}_6]^+)_i$ and $([T^{(k)}_7]^+)_i$ as
  well as in $([T^{(k-1)}_6]^+)_i$ and $([T^{(k-1)}_7]^+)_i$. At the
  same time, only two of these four terms are nonzero for any $k \in
  \KK \setminus {\horizon}$ and any $i \in \until{n}$. Using these
  facts one can obtain the bound $\inorm{\zeta_{1,a}-\zeta_{2,a}} \le
  \frac{3}{\eps}$.  Finally, the proof concludes noting
  \begin{align*}
    \inorm{& \zeta_1 - \zeta_2} = \inorm{\zeta_{1,a} + \zeta_{1,b} -
      \zeta_{2,a}- \zeta_{2,c}}
    \\
    & \qquad \le \inorm{\zeta_{1,a}-\zeta_{2,a}} + \inorm{\zeta_{1,b}}
    + \inorm{\zeta_{2,c}} = \frac{\horizon+4}{\eps}. \quad \IEEEQED
   \end{align*}
   \renewcommand{\IEEEQED}{}
\end{IEEEproof}
\renewcommand{\IEEEQED}{\IEEEQEDclosed}

 \emph{B. The \daclgg coordination algorithm:
}Here, we present our distributed algorithm and establish its
asymptotic convergence to the set of solutions of the \DEDS problem
starting from any initial condition. Our design combines ideas of
Laplacian-gradient dynamics~\cite{AC-JC:15-tcns} and dynamic average
consensus~\cite{SSK-JC-SM:15-ijrnc}.
Consider the set-valued dynamics,
\begin{subequations}\label{eq:dac-lap}
  \begin{align}
    \dot I & \in -(\eye_\horizon \otimes \Lap) \partial_I f^\eps(I,S)
    + \nu_1 z,
    \label{eq:dac-lap-1}
    \\
    \dot S & \in -\partial_S f^\eps(I,S), \label{eq:dac-lap-2}
    \\
    \dot z & = -\alpha z - \beta (\eye_\horizon \otimes \Lap) z - v +
    \nu_2(L \otimes e_r + \tilde{l} - I),\label{eq:dac-lap-3}
    \\
    \dot v & = \alpha \beta (\eye_\horizon \otimes \Lap)
    z,\label{eq:dac-lap-4}
  \end{align}
\end{subequations}
where $\alpha, \beta, \nu_2, \nu_2 \in \realpositive$ are design
parameters and $e_r \in \real^n$ is the unit vector along the $r$-th
coordinate. This dynamics is an interconnected system with two parts:
the $(I,S)$-component seeks to adjust the injection levels to satisfy
the load profile and search for the optimizers of the \DEDS problem
while the $(z,v)$-component corresponds to the dynamic average
consensus part, with $z^{(k)}_i$ aiming to track the difference
between the load $l^{(k)} = L^{(k)} + \ones_n^\top
\tilde{l}_i^{(k)}$.
Our terminology \daclgg dynamics to refer to~\eqref{eq:dac-lap} is
motivated by this ``dynamic average consensus in $(z,v)$+ Laplacian
gradient in $I$ + gradient in~$S$'' structure.  For convenience, we
denote~\eqref{eq:dac-lap} by
$\setmap{\dacLapgg}{\real^{4n\horizon}}{\real^{4n\horizon}}$.  Note
$\Eq{\dacLapgg} = \FFDEDSo$ and since $\partial_I f^\eps$ and
$\partial_S f^\eps$ are locally bounded, upper semicontinuous and take
nonempty convex compact values, the solutions of $\dacLapgg$ exist
starting from any initial condition (cf. Section~\ref{se:Prelim}).

\begin{remark}\longthmtitle{Distributed implementation of the \daclgg
    dynamics}\label{re:dist-imp}
  Writing the $(z,v)$ dynamics componentwise, one can see that for
  each $i$ and each $k$, the values $(\dot z_i^{(k)},\dot v_i^{(k)})$
  can be computed using the state variables $(z_i^{(k)},
  \{z_j^{(k)}\}_{j \in \nout(i)}, v_i^{(k)}, I_i^{(k)})$ only.
  Hence,~\eqref{eq:dac-lap-3} and~\eqref{eq:dac-lap-4} can be
  implemented in a distributed manner where each unit only
  requires information from its out-neighbors.  Subsequently, $f^\eps$
  can be written in the separable form
  \begin{align*}
    f^\eps(I,S) = \sum_{i=1}^n f^\eps_i(I^{(1)}_i, \dots,
    I^{(\horizon)}_i, S^{(1)}_i, \dots, S^{(\horizon)}_i).
  \end{align*}
  Thus, if $\zeta_1 \in \partial_I f^\eps(I,S)$ and $\zeta_2
  \in \partial_S f^\eps(I,S)$ then, for all $k \in \KK$,
  $(\zeta_1)_i^{(k)},(\zeta_2)_i^{(k)} \in \real$ only depend on the
  state of unit $i$, i.e., $(I^{(1)}_i, \dots, I^{(\horizon)}_i,
  S^{(1)}_i, \dots, S^{(\horizon)}_i)$ and are computable
  by~$i$. Hence, the $S$-dynamics can implemented by the DERs
  using their own state and to execute the $I$-dynamics, each $i$
  needs information from its out-neighbors. 
  \oprocend
\end{remark}

We next address the convergence analysis of~\eqref{eq:dac-lap}. For
convenience, let $\Hglobal = \real^{n\horizon} \times
\real^{n\horizon} \times \real^{n\horizon} \times (\HH_0)^\horizon $
and $\Homega = \prod_{k=1}^\horizon \HH_{l^{(k)}} \times
\real^{n\horizon} \times (\HH_0)^\horizon \times (\HH_0)^\horizon$.

\begin{theorem}\longthmtitle{Convergence of the \daclgg dynamics to
    the solutions of the \DEDS problem}\label{th:convergence}
  Let $\FFDEDSo$ satisfy Assumption~\ref{as:regular}, $\eps$
  satisfy~\eqref{eq:eps-bound}, and $\alpha,
  \beta, \nu_1, \nu_2 > 0$ satisfy
  \begin{equation}\label{eq:alpha-beta-cond-n}
    \frac{\nu_1}{\beta \nu_2 \lambda_2(\Lap + \Lap^\top)} +
    \frac{\nu_2^2 \lambda_{\max}(\Lap^\top \Lap)}{2 \alpha} <
    \lambda_2(\Lap + \Lap^\top).
  \end{equation}
  Then, any trajectory
  of~\eqref{eq:dac-lap} starting in $\Hglobal$
  converges to $\FFa $ $=\setdef{(I,S,z,v) \in \FFDEDSo
    \times \{0\} \times \real^{n\horizon}} {v= \nu_2 (l \otimes e_r -
    I)}$.
\end{theorem}
\begin{IEEEproof}
  Our first step is to show that the $\omega$-limit set of any
  trajectory of~\eqref{eq:dac-lap} with initial condition
  $(I_0,S_0,z_0,v_0) \in \Hglobal$ is contained in $\Homega$. To this
  end, write~\eqref{eq:dac-lap-4} as
  \begin{align*}
    \dot v^{(k)} = \alpha \beta \Lap z^{(k)} \quad \text{ for all  } k
    \in \KK.
  \end{align*}
  Note that $\ones_n^\top \dot v^{(k)} = \alpha \beta \ones_n^\top \Lap
  z^{(k)} = 0$ for all $k \in \KK$ because $\GG$ is weight-balanced.
  Therefore, the initial condition $v_0 \in (\HH_0)^\horizon$ implies
  that $v(t) \in (\HH_0)^\horizon$ for all $t \ge 0$ along any
  trajectory of~\eqref{eq:dac-lap} starting at $(I_0,S_0,z_0,v_0)$.
  Now, if $\zeta \in \partial_I f^\eps(I,S)$ then,
  from~\eqref{eq:dac-lap-1} and~\eqref{eq:dac-lap-3}, we get for any $k
  \in \KK$
  \begin{align*}
    \dot I^{(k)} & = - \Lap \zeta^{(k)} + \nu_1 z^{(k)},
    \\
    \dot z^{(k)} & = - \alpha z^{(k)} - \beta \Lap z^{(k)} - v^{(k)} +
    \nu_2 (l^{(k)} e_r - I^{(k)}).
  \end{align*}
  Let $\xi_k = \ones_n^\top I^{(k)} - l^{(k)}$. Then, from the above
  equations we get $\dot \xi_k = \ones_n^\top \dot I^{(k)} = \nu_1
  \ones_n^\top z^{(k)}$. Further, we have
  \begin{align*}
    \ddot \xi_k &  = \nu_1 \ones_n^\top \dot z^{(k)} = - \alpha \nu_1
    \ones_n^\top z^{(k)} + \nu_1 \nu_2 (l^{(k)} - \ones^\top I^{(k)})
    \\
    & = - \alpha \dot \xi_k - \nu_1 \nu_2 \xi_k,
  \end{align*}
  forming a second-order linear system for $\xi_k$. The LaSalle
  Invariance Principle~\cite{HKK:02} with the function $\nu_1 \nu_2
  \norm{\xi_k}^2 + \norm{\dot \xi_k}^2$ implies that as $t \to \infty$
  we have $(\xi_k(t); \dot \xi_k(t)) \to 0$ and so $\ones_n^\top
  I^{(k)}(t) \to l^{(k)}$ and $\ones_n^\top z^{(k)}(t) \to 0$ as $t \to
  \infty$.

  Next, proceeding to the convergence analysis, consider the change of
  coordinates  $\map{D}{\real^{4n\horizon}}{\real^{4n\horizon}}$ defined by
  \begin{align*}
    (I,S,\omega_1,\omega_2) & = D(I,S,z,v) \\ & = (I,S,z,v +\alpha z-
    \nu_2 (l \otimes e_r -I)) .
  \end{align*}
  In these coordinates, the set-valued map~\eqref{eq:dac-lap}
  takes the form
  \begin{align}
    \dacLapgg &(I,S,\omega_1,\omega_2) = \setdef{( -(\eye_\horizon \otimes \Lap
      )\zeta_1 + \nu_1 \omega_1, -\zeta_2, \notag
      \\
      &  -\beta(\eye_\horizon \otimes \Lap) \omega_1 - \omega_2,
      \label{eq:rewrite-edp2-2}
      \\
      & \nu_1 \nu_2 \omega_1 -\alpha \omega_2 
      - \nu_2 (\eye_\horizon \otimes \Lap) \zeta_1) \in
      \real^{4n\horizon}}{ \notag
      \\
      & \zeta_1 \in \partial_I f^\eps (I,S), \zeta_2 \in \partial_S
      f^\eps (I,S)} . \notag
  \end{align}
  This transformation helps in identifying the LaSalle-type function
  for the dynamics.  We now focus on proving that, in the new
  coordinates, the trajectories of~\eqref{eq:dac-lap} converge
  to 
  \begin{align*}
    \bFFa & = D(\FFa) = \FFDEDSo \times \{0\} \times \{0\}.
  \end{align*}
  Note that $D(\Homega) = \Homega$ and so, from the property of the
  $\omega$-limit set of trajectories above, we get that $t \mapsto
  (I(t),S(t),\omega_1(t),\omega_2(t))$ starting in $D(\Hglobal)$
  belongs to~$\Homega$.  Next, we show the hypotheses of
  Proposition~\ref{pr:refined-lasalle-nonsmooth} are satisfied, where
  $\Homega$ plays the role of $\SS \subset \real^{4n\horizon}$
  and 
  $\map{V}{\real^{4n\horizon}}{\realnonnegative}$,
  \begin{align*}
    V(I,S,\omega_1,\omega_2) = f^\eps (I,S) + \tfrac{1}{2}(\nu_1 \nu_2
    \norm{\omega_1}^2 + \norm{\omega_2}^2).
  \end{align*} 
  plays the role of~$W$, resp.
  Let $(I,S,\omega_1,\omega_2) \in \Homega$ then any element of
  $\SetLie_{\dacLapgg}V(I,S,\omega_1,\omega_2)$ can be written as
  \begin{align}\label{eq:liederivative-element}
    & - \zeta_1^\top (\eye_\horizon \otimes \Lap) \zeta_1 + \nu_1
    \zeta_1^\top \omega_1 - \norm{\zeta_2}^2 
    - \beta \nu_1 \nu_2 \omega_1^\top (\eye_\horizon \otimes \Lap)
    \omega_1 \notag 
    \\
    & - \alpha \norm{\omega_2}^2 
    - \nu_2 \omega_2^\top (\eye_\horizon \otimes L) \zeta_1,
  \end{align}
  where $\zeta_1 \in \partial_I f^\eps (I,S)$ and $\zeta_2 \in
  \partial_S f^\eps (I,S)$. Since the digraph $\GG$ is strongly
  connected and weight-balanced, we use~\eqref{eq:LapBound} and
  $\ones_{n\horizon}^\top \omega_1 = 0$ to bound the above expression
  as
  \begin{align*}
    & - \tfrac{1}{2}\lambda_2(\Lap + \Lap^\top) \norm{\eta}^2 + \nu_1
    \eta^\top \omega_1 - \norm{\zeta_2}^2
    \\
    & - \tfrac{1}{2} \beta \nu_1 \nu_2 \lambda_2(\Lap + \Lap^\top)
    \norm{\omega_1}^2 \notag - \alpha \norm{\omega_2}^2 - \nu_2
    \omega_2^\top (\eye_\horizon \otimes \Lap) \eta
    \\
    & = \gamma^\top M \gamma - \norm{\zeta_2}^2 ,
  \end{align*}
  where $\eta = (\eta^{(1)}; \dots; \eta^{(\horizon)})$ with
  $\eta^{(k)} = \zeta^{(k)} - \tfrac{1}{n} (\ones_n^\top
  \zeta^{(k)})\ones_n$, the vector $\gamma = (\eta; \omega_1;
  \omega_2)$, and the matrix
  \begin{align*} 
    M = \begin{bmatrix} -\tfrac{1}{2} \lambda_2(\Lap + \Lap^\top)
      \eye_{n\horizon} & B^\top \\ B & C
    \end{bmatrix}, 
  \end{align*}
  with $B^\top = \begin{bmatrix}
    \tfrac{1}{2} \nu_1 \eye_{n\horizon} & - \tfrac{1}{2} \nu_2
    (\eye_{\horizon} \otimes \Lap)^\top
  \end{bmatrix}$, and
  \begin{align*}
    C =\begin{bmatrix} -\tfrac{1}{2} \beta \nu_1 \nu_2 \lambda_2(\Lap
      + \Lap^\top) \eye_{n\horizon} & 0
      \\
      0 & -\alpha \eye_{n\horizon}
    \end{bmatrix}.
  \end{align*}
  Resorting to the Schur complement~\cite{SB-LV:09}, $M \in
  \real^{3n\horizon \times 3n\horizon}$ is neg. definite if $
  -\tfrac{1}{2} \lambda_2(\Lap + \Lap^\top) \eye_{n\horizon} - B^\top
  C^{-1} B$, that equals
  \begin{align*}
     -\tfrac{1}{2} \lambda_2(\Lap + \Lap^\top)\eye_{n\horizon} +
    \tfrac{\nu_1}{2 \beta \nu_2 \lambda_2(\Lap + \Lap^\top)}
    \eye_{n\horizon}+ \tfrac{\nu_2^2}{4\alpha} (\eye_{\horizon} \otimes
    \Lap)^\top (\eye_{\horizon} \otimes \Lap) ,
  \end{align*}
  is negative definite, which follows
  from~\eqref{eq:alpha-beta-cond-n}. Hence, for any
  $(I,S,\omega_1,\omega_2) \in \Homega$, we have $\max
  \SetLie_{\dacLapgg} V(I,S,\omega_1,\omega_2) \le 0$ and also $0 \in
  \SetLie_{\dacLapgg} V(I,S,\omega_1,\omega_2)$ iff $\eta = \zeta_2 =
  \omega_1 = \omega_2 =0$, which means $\zeta^{(k)} \in
  \spn\{\ones_n\}$ for each $k \in \KK$.  Consequently, using the
  characterization of optimizers in~\eqref{eq:solution-set}, we deduce
  that $(I,S)$ is a solution of~\eqref{eq:deds-m} and so,
  $(I,S,\omega_1,\omega_2) \in \bFFa$. Since, $\bFFa$ belongs to a
  level set of~$V$, we conclude that
  Proposition~\ref{pr:refined-lasalle-nonsmooth}\ref{as:refined-lasalle-1}
  holds. Further, using~\cite[Lemma A.1]{AC-JC:14-auto} one can show
  that
  Proposition~\ref{pr:refined-lasalle-nonsmooth}\ref{as:refined-lasalle-2}
  holds too (we omit the details due to space constraints).

  To apply Proposition~\ref{pr:refined-lasalle-nonsmooth}, it remains
  to show that the trajectories starting from $D(\Hglobal)$ are
  bounded.  We reason by contradiction. Assume there exists $t \mapsto
  (I(t), S(t), \omega_1(t),\omega_2(t))$, with
  $(I(0),S(0),\omega_1(0),\omega_2(0)) \in D(\Hglobal)$, of
  $\dacLapgg$ such that $\norm{(I(t),S(t),\omega_1(t),\omega_2(t)} \to
  \infty$.  Since $V$ is radially unbounded, this implies
  $V(I(t),S(t),\omega_1(t),\omega_2(t)) \to \infty$.  Also, as
  established above, we know $\ones_n^\top I^{(k)}(t) \to l^{(k)}$ and
  $\ones_n^\top \omega_1^{(k)}(t) \to 0$ for each $k \in \KK$.  Thus,
  there exist times $\{t_m\}_{m=1}^{\infty}$ with $t_m \to \infty$
  such that for all $m \in \integerspositive$,
  \begin{align}\label{eq:xi-bound}
    \abs{\ones_n^\top \omega_1^{(k)}(t_m)}  < {1}/{m} \text{ for all
    }  k \in & \KK
    , 
    \\
    \max \SetLie_{\dacLapgg} V(I(t_m),S(t_m),\omega_1(t_m),
    \omega_2(t_m)) & >
    0. \notag 
  \end{align}
  The second inequality implies the existence of $\{\zeta_{1,m}\}_{m =
    1}^{\infty}$ and $\{\zeta_{2,m}\}_{m=1}^{\infty}$ with
  $(\zeta_{1,m},\zeta_{2,m}) \in (\partial_I
  f^{\eps}(I(t_m),S(t_m)), \partial_S f^{\eps}(I(t_m),S(t_m)))$, such
  that
  \begin{align*}
    - \zeta_{1,m}^\top & (\eye_\horizon \otimes \Lap) \zeta_{1,m} + \nu_1 
    \zeta_{1,m}^\top \omega_1(t_m)  - \norm{\zeta_{2,m}}^2  \notag
    \\
    & - \beta \nu_1 \nu_2 \omega_1(t_m)^\top (\eye_\horizon \otimes \Lap)
    \omega_1(t_m) - \alpha \norm{\omega_2(t_m)}^2 \notag
    \\
    & - \nu_2 \omega_2(t_m)^\top (\eye_\horizon \otimes \Lap) \zeta_{1,m} > 0,
  \end{align*}
  for all $m \in \integerspositive$, where we
  have used~\eqref{eq:liederivative-element} to write an element of 
  $\SetLie_{\dacLapgg} V (I,S,\omega_1,\omega_2)$. Letting
  $\eta_m^{(k)} = \zeta_{1,m}^{(k)} - \tfrac{1}{n}(\ones_n^\top
  \zeta_{1,m}^{(k)}) \ones_n$, using~\eqref{eq:LapBound}, and using the
  relation $\norm{\omega_1^{(k)}(t_m) - \tfrac{1}{n} (\ones_n^\top
      \omega_1^{(k)}(t_m)) \ones_n}^2 = \norm{\omega_1^{(k)}(t_m)}^2
    - \tfrac{1}{n} (\ones_n^\top \omega_1^{(k)}(t_m))^2$, the above
    inequality can be rewritten as 
  \begin{multline}\label{eq:auxxx}
    \gamma_m^\top M \gamma_m + \tfrac{1}{n} \nu_1 \sum_{k\in \KK}
    (\ones_n^\top \zeta_{1,m}^{(k)})(\ones_n^\top \omega_1^{(k)}(t_m))
    - \norm{\zeta_{2,m}}^2
    \\
    + \tfrac{\beta \nu_1 \nu_2}{2 n } \lambda_2(\Lap + \Lap^\top)
    \sum_{k \in \KK} (\ones_n^\top \omega_1^{(k)}(t_m))^2 > 0,
  \end{multline}
  with $\gamma_m = (\eta_{m}; \omega_1(t_m); \omega_2(t_m))$.
  Using~\eqref{eq:xi-bound} on~\eqref{eq:auxxx},
  \begin{multline}\label{eq:lie-bound}
    \gamma_m^\top M \gamma_m - \norm{\zeta_{2,m}}^2 + \tfrac{\nu_1}{nm}
    \sum_{k \in \KK} \abs{\ones_n^\top \zeta_{1,m}^{(k)}}
    \\
    + \tfrac{\beta \nu_1 \nu_2 \horizon}{2 n m^2} \lambda_2(\Lap +
    \Lap^\top) > 0
  \end{multline}
  for all $m \in \integerspositive$. 
  Next, we consider two cases, depending on whether the sequence
  $\{(I(t_m),S(t_m))\}_{m=1}^\infty$ is (a) bounded or (b)
  unbounded.  
  In case (a), $\{(\omega_1(t_m),\omega_2(t_m))\}_{m=1}^\infty$ must
  be unbounded. Since $M$ is negative definite, we have $\gamma_m^\top
  M \gamma_m \le \lambda_{\max}(M) \norm{(\omega_1(t_m),
    \omega_2(t_m))}^2$.  Thus, by~\eqref{eq:lie-bound}
    \begin{align}\label{eq:omega-inf-bound}
    \lambda_{\max} (M) \norm{(\omega_1(t_m),& \omega_2(t_m))}^2  +
    \tfrac{\nu_1}{nm} \sum_{k \in \KK} \abs{\ones_n^\top
      \zeta_{1,m}^{(k)}} \notag
    \\
    & + \tfrac{\beta \nu_1 \nu_2 \horizon}{2 n m^2} \lambda_2(\Lap +
    \Lap^\top) > 0.
  \end{align} 
  Since $\partial_I f^\eps$ is locally bounded and
  $\{(I(t_m),S(t_m))\}_{m=1}^\infty$ is bounded, we deduce
  $\{\zeta_{1,m}\}$ is bounded~\cite[Proposition 6.2.2]{JBHU-CL:93}.
  Combining these facts with $\lambda_{\max} (M) <0$ and
  $\norm{(\omega_1(t_m),\omega_2(t_m))} \to \infty$, one can find
  $\bar{m} \in \integerspositive$ such that~\eqref{eq:omega-inf-bound}
  is violated for all $m \ge \bar{m}$, a contradiction.  Now consider
  case (b) where $\{(I(t_m),S(t_m))\}_{m=1}^\infty$ is unbounded. We
  divide this case further into two, based on the sequence
  $\bigl\{\sum_{k=1}^\horizon \abs{\ones_n^\top \zeta_{1,m}^{(k)}}
  \bigr\}_{m=1}^\infty$ being bounded or not. Using $\gamma_m^\top M
  \gamma_m \le \lambda_{\max}(M)\norm{\eta_m}^2$, the
  inequality~\eqref{eq:lie-bound} implies
  \begin{multline}\label{eq:state-bound}
    \lambda_{\max} (M) \norm{\eta_m}^2 - \norm{\zeta_{2,m}}^2 +
    \frac{\nu_1}{nm} \sum_{k=1}^\horizon \abs{\ones_n^\top
      \zeta_{1,m}^{(k)}}
    \\
    + \frac{\beta \nu_1 \nu_2 \horizon}{2 n m^2} \lambda_2(\Lap +
    \Lap^\top) > 0.
  \end{multline}
  Consider the case when $\bigl\{\sum_{k=1}^\horizon \abs{\ones_n^\top
    \zeta_{1,m}^{(k)}} \bigr\}_{m=1}^\infty$ is unbounded. Partition
  $\KK$ into disjoint sets $\KK_u$ and $\KK_b$ such that
  $\abs{\ones_n^\top \zeta_{1,m}^{(k)}} \to \infty$ for all $k \in
  \KK_u$ and $\bigl\{\abs{\ones_n^\top \zeta_{1,m}^{(k)}}\bigr\}_{m=1}
  ^\infty$ is uniformly bounded for all $k \in \KK_b$.  For
  convenience, rewrite~\eqref{eq:state-bound} as $\sum_{k =
    1}^\horizon U_{k,m} + \frac{Z_1}{m} > 0$, where $ Z_1 =
  \frac{\beta \nu_1 \nu_2 \horizon}{2nm} \lambda_2(\Lap + \Lap^\top)$
  and, for each $k \in \KK$,
  \begin{multline*}
    U_{k,m} = \lambda_{\max}(M) \norm{\eta_m^{(k)}}^2 -
    \norm{\zeta_{2,m}^{(k)}}^2 + \frac{\nu_1}{nm} \abs{\ones_n^\top
      \zeta_{1,m}^{(k)}}.
  \end{multline*}
  By definition of~$\KK_b$, there exists $Z_2 > 0$ with $\sum_{k \in
    \KK_b} U_{k,m} \le \frac{Z_2}{m}$. Hence,
  if~\eqref{eq:state-bound} holds for all $m \in \integerspositive$,
  then so is
  \begin{align*}
    \sum_{k \in \KK_u} U_{k,m} + \frac{Z_1 + Z_2}{m} > 0.
  \end{align*}
  Next we show that for each $k \in \KK_u$ there exists $m_k \in
  \integerspositive$ such that $U_{k,m} + \frac{Z_1 + Z_2}{m} < 0$ for
  all $m \ge m_k$. This will lead to the desired contradiction.
  Assume without loss of generality that $\ones_n^\top
  \zeta_{1,m}^{(k)} \to \infty$ (reasoning for the case when the
  sequence approaches negative infinity follows analogously). Then,
  for
  \begin{align*}
    \lambda_{\max}(M) \norm{\eta_m^{(k)}}^2 - \norm{\zeta_{2,m}}^2 +
    \frac{\nu_1}{nm} \abs{\ones_n^\top \zeta_{1,m}^{(k)}} + \frac{Z_1 +
    Z_2}{m} > 0,
  \end{align*}
  for all $m \in \integerspositive$, we require $(\zeta_{1,m}^{(k)})_i
  \to \infty$ for all $i \in \until{n}$.
  Indeed, otherwise, recalling that $\eta_m^{(k)} = \zeta_{1,m}^{(k)}
  - \frac{1}{n}(\ones_n^\top \zeta_{1,m}^{(k)}) \ones_n$, it can be
  shown that there exist an $\bar{m}$ such that 
  \begin{align*}
    \lambda_{\max} \norm{\eta_m^{(k)}}^2 < \frac{\nu_1}{nm} 
    \abs{\ones_n^\top \zeta_{1,m}^{(k)}} + \frac{Z_1 +
    Z_2}{m} \, \,  \text{ for all } m \ge \bar{m}.
  \end{align*}
  Note that from Lemma~\ref{le:size-gen-gradient} we have
  $\inorm{\zeta_{1,m}^{(k)} - \zeta_{2,m}^{(k)}} \le
  \frac{\horizon+4}{\eps}$ which further implies that
  $(\zeta_{2,m}^{(k)})_i \to \infty$ for all $i \in \until{n}$.  With
  these facts in place, we write
  \begin{align*}
    U_{k,m} + \frac{Z_1 + Z_2 }{m} & < - \sum_{i=1}^n
    (\zeta_{2,m}^{(k)})_i^2 +
    \frac{\nu_1}{m} \abs{\sum_{i=1}^n (\zeta_{1,m}^{(k)})_i}
    \\
    & \qquad + \frac{Z_1 + Z_2}{m}
  \end{align*}
  and deduce that there exists an $m_k \in \integerspositive$ such
  that the right-hand side of the above expression is negative for all
  $m \ge m_k$, which is what we wanted to show.
  
  Finally, consider the case when the sequence $\bigl
  \{\sum_{k=1}^\horizon \abs{\ones_n^\top \zeta_{1,m}^{(k)}} \bigr
  \}_{m=1}^\infty$ is bounded.  For~\eqref{eq:state-bound} to be true
  for all $m \in \integerspositive$, we require $\norm{\gamma_m} \to
  0$ and $\norm{\zeta_{2,m}} \to 0$ as $m \to \infty$.  This further
  implies that $\eta_m \to 0$ and, from Assumption~\ref{as:regular},
  this is only possible if $\{(I(t_m),S(t_m))\}_{m=1}^\infty$ is
  bounded, which is a contradiction.
\end{IEEEproof}

\begin{remark}\longthmtitle{General setup for storage:
    revisited}\label{re:storage-subset-re}
  The \daclgg dynamics~\eqref{eq:dac-lap} can be modified to scenarios
  that include more general descriptions of storage capabilities, as
  in Remark~\ref{re:storage-subset}. For instance, if only a subset of
  units have storage capabilities, the only modification is to set the
  variables $\{S_i^{(k)}\}_{i \in \vertices_{g}, k\in \KK}$ to zero
  and execute~\eqref{eq:dac-lap-2} only for the variables
  $\{S_i^{(k)}\}_{i \in \vertices_{gs}, k\in \KK}$. 
  The resulting strategy converges to the solution set of the
  corresponding \DEDS problem.  \oprocend
\end{remark}

\begin{remark}\longthmtitle{Distributed selection of design
    parameters}\label{re:dist-selection}
  The implementation of the \daclgg dynamics requires the selection of
  parameters $\alpha,\beta,\nu_1,\nu_2,\eps$
  satisfying~\eqref{eq:eps-bound} and~\eqref{eq:alpha-beta-cond-n}.
  Condition~\eqref{eq:alpha-beta-cond-n} involves knowledge of
  network-wide quantities, but the units can resort to various
  distributed procedures to collectively select appropriate values.
  Regarding~\eqref{eq:eps-bound}, an upper bound on the denominator of
  the right-hand side can be computed aggregating, using consensus,
  the difference between the max and the min values that each
  DER's aggregate cost function takes in its respective
  feasibility set (neglecting load conditions).  The challenge for the
  units, however, is to estimate the parameter~$\rho$ if it is not
  known a priori.
  \oprocend
\end{remark}

\myclearpage 
\section{Simulations}\label{sec:sims}

We illustrate the application of the \daclgg dynamics to solve the
\DEDS problem for a group of $n=10$ generators with communication
defined by a directed ring with bi-directional edges
$\{(1,5),(2,6),(3,7),(4,8)\}$ (all edge weights are $1$).
The planning horizon is $\horizon=6$ and the load profile consists of
the external load $L = (1950, 1980, 2700, 2370, 1900,1850)$ and the
load at each generator $i$ for each slot $k$ given by
$\tilde{l}_i^{(k)} = 10i$. Thus, for each slot $k$, $\tilde{l}^{(k)} =
\sum_{i=1}^{10} \tilde{l}_i^{(k)} = 550$ and so, $l = (2500, 2530,
3250, 2920, 2450, 2400)$.
Generators have storage capacities determined by $C^M = 100 \ones_n$
and $C^m = S^{(0)} = 5 \ones_n$.  The cost function of each unit is
quadratic and constant across time.  Table~\ref{tb:Cost1} details the
cost function coefficients, generation limits, and ramp constraints,
which are modified from the data for $39$-bus New England
system~\cite{RDQ-CEMS-RJT:11}.

\begin{table}[htb]
  \centering
  {\small
    \begin{tabular}{ | c | c  c  c  c  c c c| } 
      \hline
      Unit & $a_i$ & $b_i$  & $c_i$ & $P_i^m$ & $P_i^M$ & $R^l_i$ & $R^u_i$
      \\
      \hline
      1     & 240 & 7.0 & 0.0070 & 0 & 1040 & 120 & 80 \\
      2     & 200 & 10.0 & 0.0095 & 0 & 646 & 90 & 50\\
      3     & 220 & 8.5 & 0.0090 & 0 & 725 & 100 & 65\\
      4     & 200 & 11.0 & 0.0090 & 0 & 652 & 90 & 50\\
      5     & 220 & 10.5 & 0.0080 & 0 & 508 & 90 & 50\\
      6     & 190 & 12.0 & 0.0075 & 0 & 687 & 90 & 50\\
      7     & 200 & 10.0 & 0.0100 & 0 & 580 & 120 & 80\\
      8     & 170 & 9.0 & 0.0090 & 0 & 564 & 90 & 50\\
      9     & 190 & 11.0 & 0.0072 & 0 & 865 & 100 & 65\\
      10     & 220 & 8.8 & 0.0080 & 0 & 1100 & 90 & 50\\
      \hline
    \end{tabular}}
  \caption{Cost coefficients $(a_i, b_i, c_i)$ and
    bounds $P_i^M$, $P_i^m, R^l_i$, $R^u_i$.  
    The cost function of  $i$ is $f_i(P_i) = a_i + b_i P_i + c_i P_i^2$.}\label{tb:Cost1}
  \vspace*{-5ex}
\end{table}

\begin{figure}[htb]
  \centering
    \subfloat[Total injection]
    {\includegraphics[height = 0.23
    \linewidth]{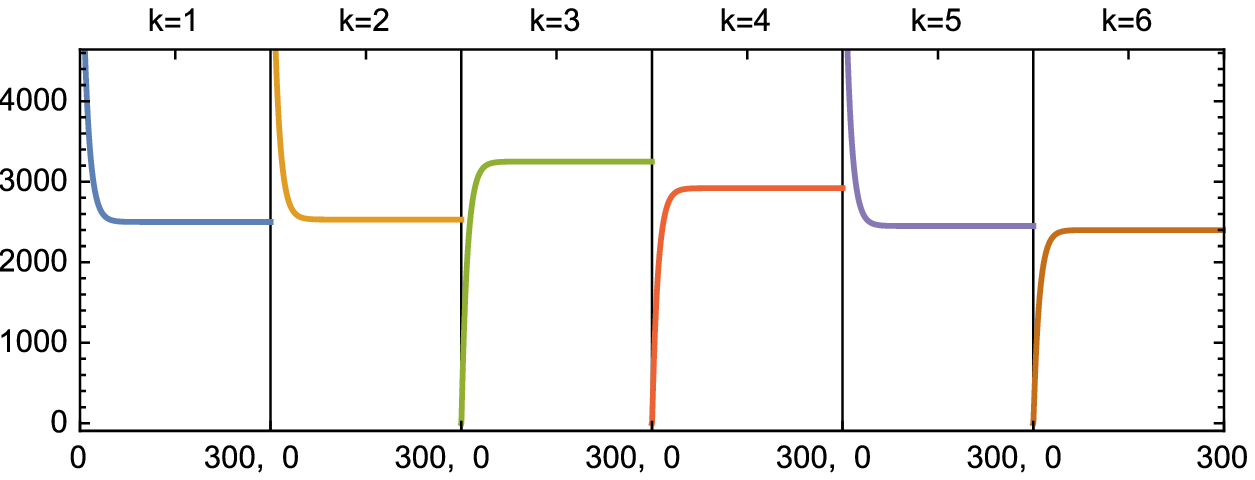}}
    \subfloat[Total cost]{\includegraphics[height =
      0.22
    \linewidth]{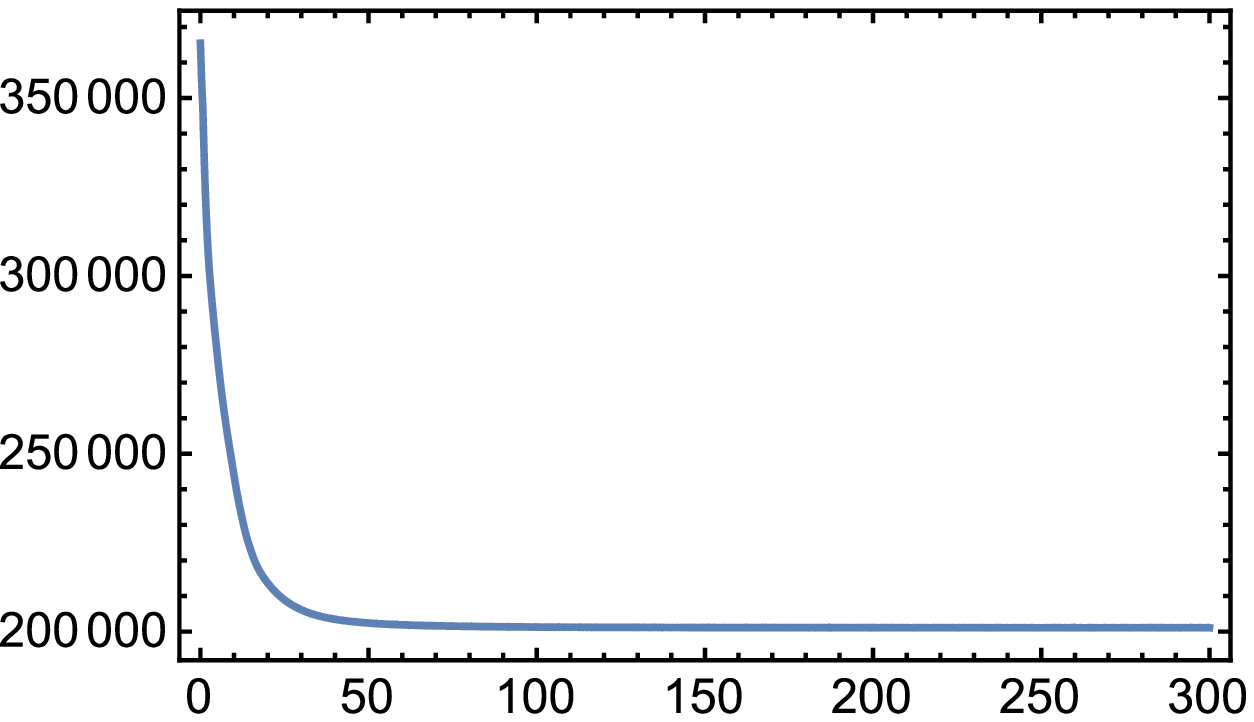}}
  \caption{Illustration of the execution of \daclgg dynamics for a
    network of $10$ generators with communication topology given by a
    directed ring among the generators with bi-directional edges
    $\{(1,5),(2,6),(3,7),(4,8)\}$ where all edge weights are $1$.
    Table~\ref{tb:Cost1} gives the box constraints, the ramp
    constraints, and the cost functions.  The load profile is $l = 
    (2500, 2530, 3250, 2920, 2450, 2400)$ and $C^M = 100 \ones_n$, $C^m
    = S^{(0)} = 5\ones_n$. Plots (a) and (b) show the time evolution
    of the total injection at each time slot and the aggregate cost
    along a trajectory of the \daclgg dynamics starting at $I(0) =
    (P^M,P^M,P^m,P^m,P^M,P^m)$, $S(0) = z(0) = v(0) =
    \zeros_{n\horizon}$. The parameters are $\eps = 0.007$, $\alpha =
    4$, $\beta = 10$, and $\nu_1 = \nu_2 = 0.65$ (which satisfy
    conditions~\eqref{eq:eps-bound}
    and~\eqref{eq:alpha-beta-cond-n}).}\label{fig:evolution}
\end{figure}

\begin{figure*}[htb]
  \centering 
   \subfloat[Power generation]{\includegraphics[width = 0.27
    \linewidth]{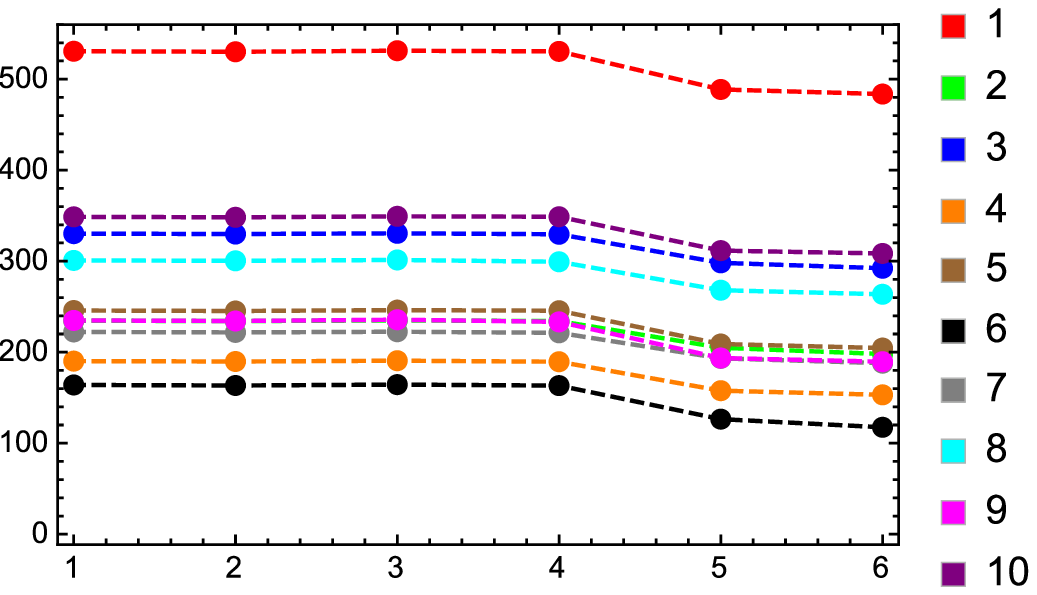}}
    \quad 
  \subfloat[Power injected into the grid]{\includegraphics[width = 0.27
    \linewidth]{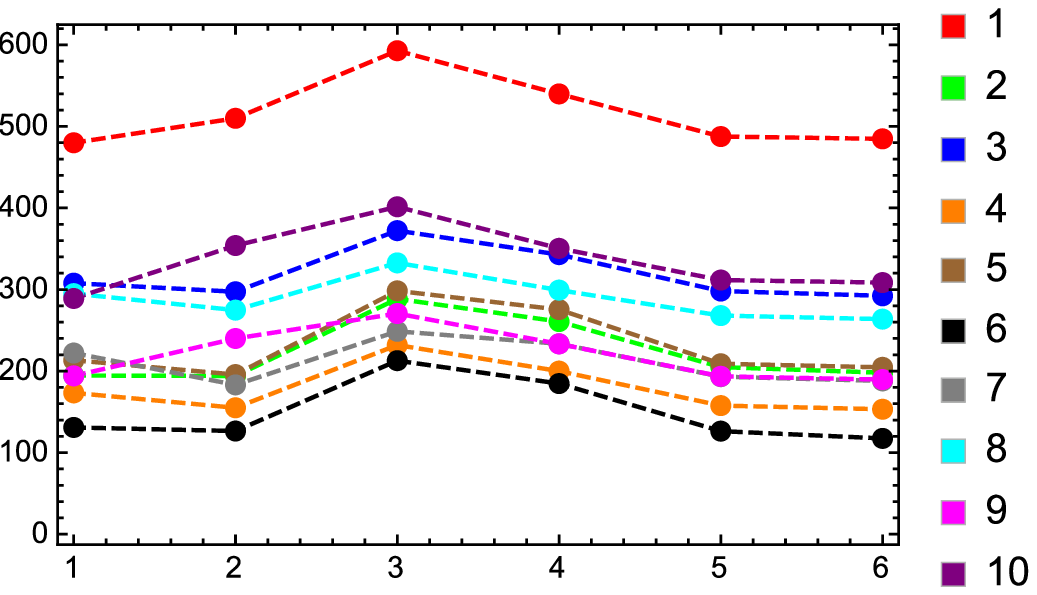}}
    \quad
  \subfloat[Power sent to storage]{\includegraphics[width = 0.27
    \linewidth]{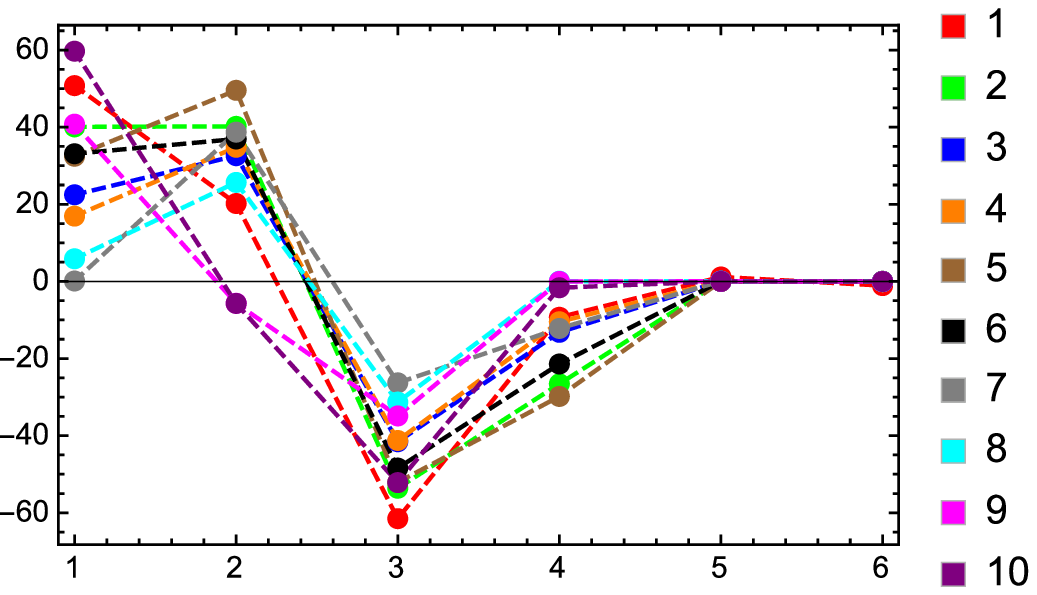}}
  \\
  \hspace{-4ex}
  \subfloat[Total power generated]{\includegraphics[width = 0.242
    \linewidth]{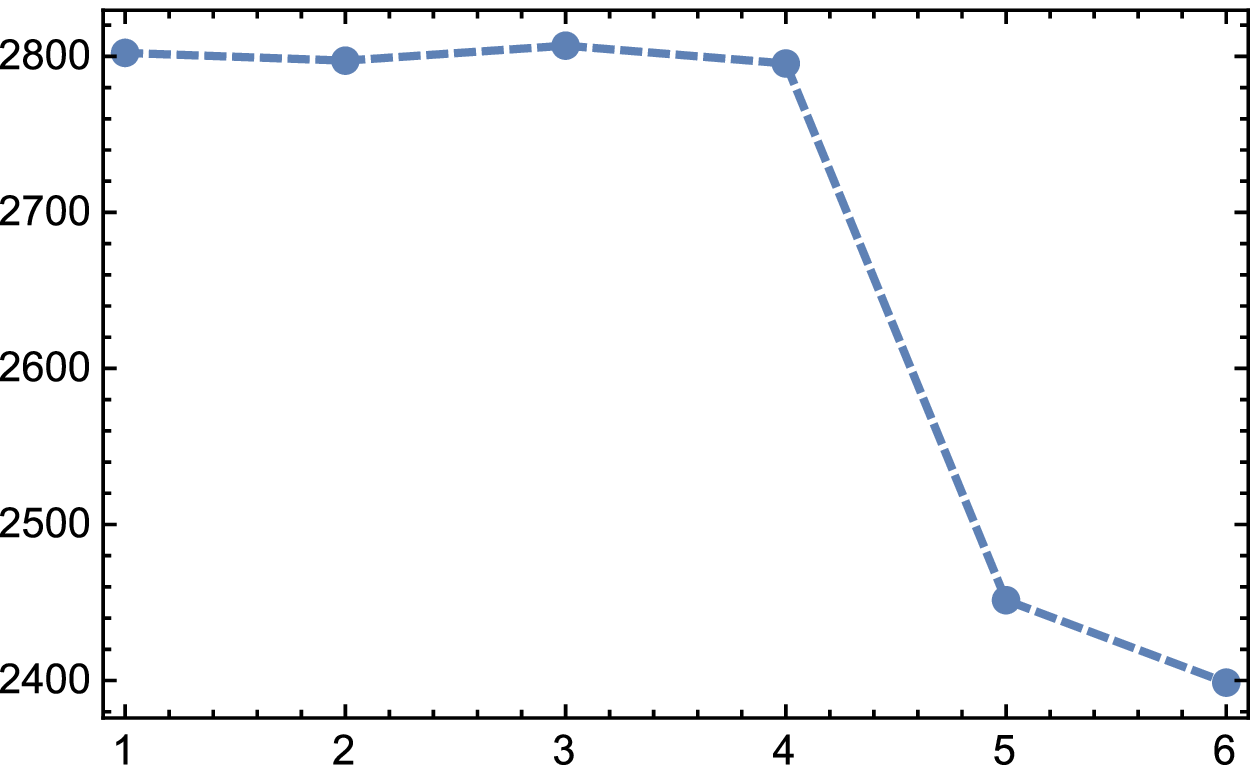}}
    \qquad \hspace{0.3ex}
  \subfloat[Total power sent to storage]{\includegraphics[width = 0.242
    \linewidth]{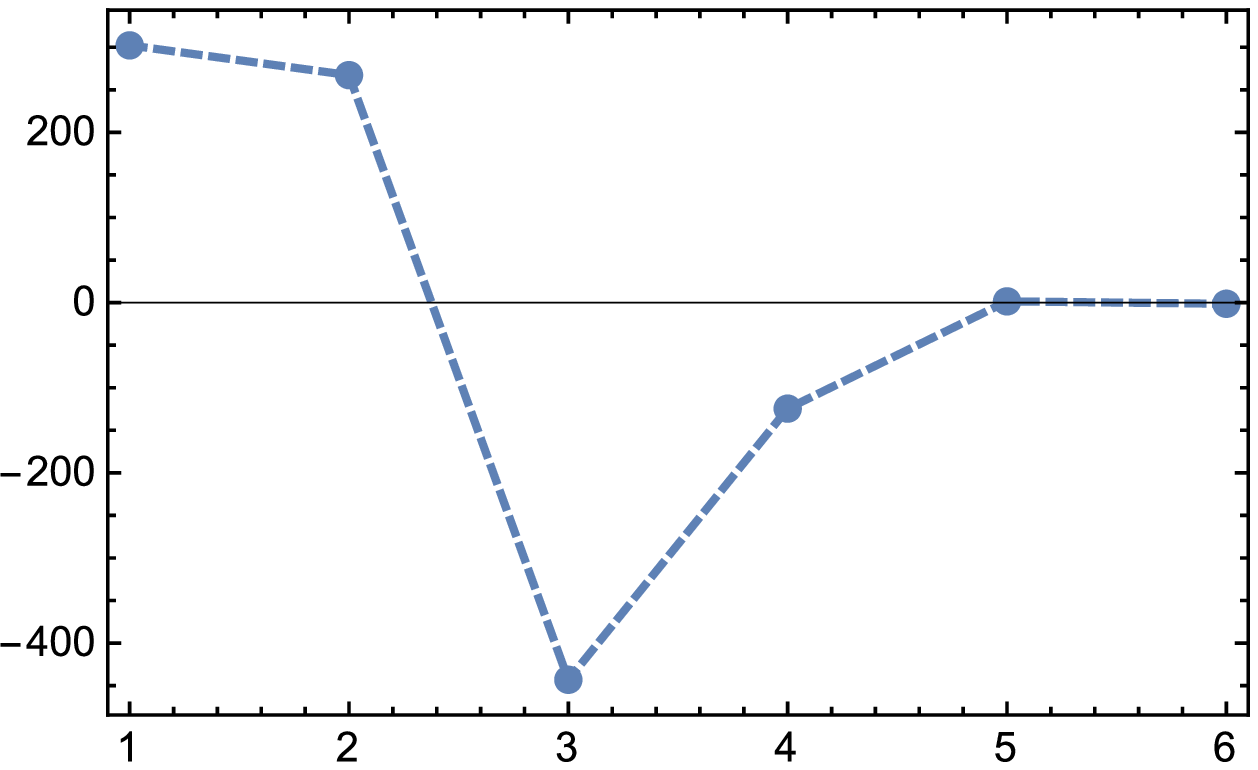}}
    \qquad \hspace{0.3ex} 
  \subfloat[Aggregate state of stored power]{\includegraphics[width =
    0.242
    \linewidth]{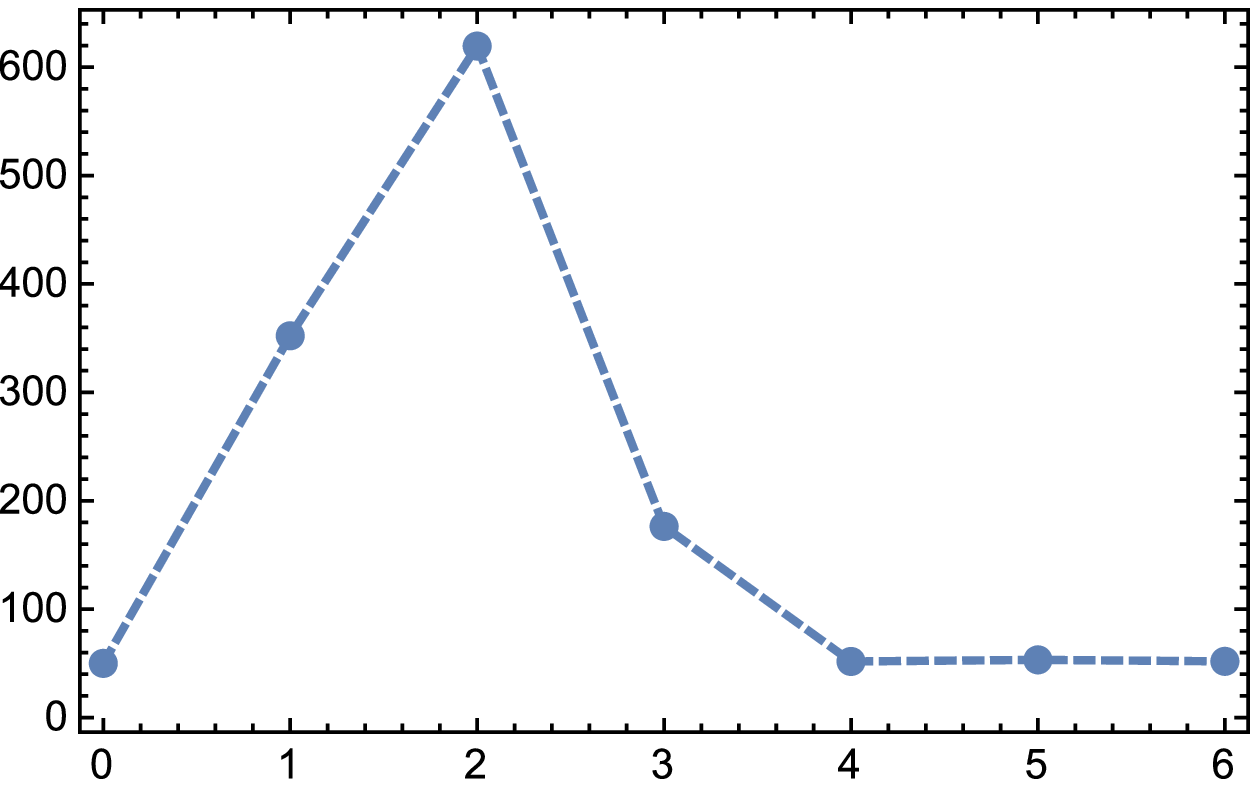}}
  \caption{Plots (a) to (f) illustrate the solution obtained in
    Figure~\ref{fig:evolution}. Plots (b) and (c) show the power
    injected and power sent to storage across the time horizon, with
    unique colors for each generator. These values add up to the total
    generation in (a). The collective behavior is represented in
    (d)-(f), where we plot the total power generated, the total power
    sent to storage, and the aggregate of the power stored in the
    storage units, respectively. The profile of total injection is the
    same as that of load profile.  Since the time-independent cost is
    quadratic with positive coefficients and the storage capacity is
    large enough, one can show that the optimal strategy is to produce
    the same power, i.e., $\tfrac{1}{5}\sum_{k=1}^5 l^{(k)}$, at each
    time slot $k = 1, \dots, 5$, as seen in (a) and (d).  The initial
    excess generation (due to the lower required load) at slots $k =
    1,2$ is stored and used in slots $k=3,4,5,6$, as indicated in (e)
    and (f).}\label{fig:conv}
\end{figure*}

Figure~\ref{fig:evolution} illustrates the evolution of the total
power injected at each time slot and the total cost incurred by the
network, respectively. As established in Theorem~\ref{th:convergence}
and shown in Figure~\ref{fig:conv}, the total injection asymptotically
converges to the load profile $l$, the total aggregate cost converges
to the minimum $201092$ and the converged solution
satisfies~\eqref{eq:box-cons}-\eqref{eq:ramp-cons}.

\section{Conclusions}

We have studied the \DEDS problem for a group of generators with
storage capabilities that communicate over a strongly connected,
weight-balanced digraph.  Using exact penalty functions, we have
provided an alternative problem formulation, upon which we have built
to design the distributed \daclgg dynamics.  This dynamics provably
converges to the set of solutions of the problem from any initial
condition.  For future work, we plan to extend the scope of our
formulation to include power flow equations, constraints on the power
lines, various losses, and stochasticity of the available data (loads,
costs, and generator availability).
We also intend to explore the use of our dynamics as a building block
in solving grid control problems across different time scales (e.g.,
implementations at long time scales on high-inertia generators and at
short time scales on low-inertia generators 
in the face of highly-varying demand) and hierarchical levels (e.g.,
in multi-layer architectures where aggregators at one layer coordinate
their response to a request for power production, and feed their
decisions as load requirements to the devices in lower layers).

\bibliographystyle{ieeetr}%
\bibliography{alias,Main,Main-add,JC}

\end{document}